\documentclass{aims}
\usepackage{amsmath,amssymb,amsfonts}
\usepackage{paralist}
\usepackage{mathrsfs}
\usepackage{graphicx}
\usepackage{cite}

\usepackage[bookmarks=false,%
pdfstartview=FitH,linkbordercolor={0.5 1 1},%
citebordercolor={0.5 1 0.5},unicode,%
pagebackref,hyperindex]{hyperref}

\def\currentvolume{X}
 \def\currentissue{X}
  \def\currentyear{200X}
   \def\currentmonth{XX}
    \def\ppages{X--XX}

\newtheorem{theorem}{Theorem}[section]
\newtheorem{corollary}{Corollary}
\newtheorem*{main}{Main Theorem}
\newtheorem{lemma}[theorem]{Lemma}

\newtheorem*{problem}{Problem}
\theoremstyle{definition}

\newtheorem{example}{Example}

\newcommand{\ecc}{\mathop{\mathrm{ecc}}}
\newcommand{\setA}{\mathscr{A}}
\newcommand{\tr}{\mathop{\mathrm{tr}}}
\newcommand{\conv}{\mathop{\mathrm{conv}}}

\title[Iterative building of Barabanov norms]%
{Iterative building of Barabanov norms and
computation of the joint spectral radius for matrix sets}

\author[Victor Kozyakin]{}


\email{kozyakin@iitp.ru}

\def\publname{}
\def\volinfo{}
\def\pageinfo{}

\begin{document}
\maketitle

\centerline{\scshape Victor Kozyakin}
\medskip
{\footnotesize
 \centerline{Institute for Information Transmission Problems}
 \centerline{Russian Academy of Sciences}
  \centerline{Bolshoj Karetny lane 19, Moscow 127994 GSP-4, Russia}
} 

\bigskip

\begin{abstract}
The problem of construction of Barabanov norms for analysis of
properties of the joint (generalized) spectral radius  of
matrix sets has been discussed in a number of publications. In
\cite{Koz:CDC05:e,Koz:INFOPROC06:e} the method of Barabanov
norms was the key instrument in disproving the Lagarias-Wang
Finiteness Conjecture. The related constructions were
essentially based on the study of the geometrical properties of
the unit balls of some specific Barabanov norms. In this
context the situation when one fails to find among current
publications any detailed analysis of the geometrical
properties of the unit balls of Barabanov norms looks a bit
paradoxical. Partially this is explained by the fact that
Barabanov norms are defined nonconstructively, by an implicit
procedure. So, even in simplest cases it is very difficult to
visualize the shape of their unit balls. The present work may
be treated as the first step to make up this deficiency. In the
paper an iteration procedure is considered that allows to build
numerically Barabanov norms for the irreducible matrix sets and
simultaneously to compute the joint spectral radius of these
sets.
\end{abstract}

\section{Introduction}\label{S-intro}

Let $\setA=\{A_{1},\ldots,A_{r}\}$ be a set of real $m\times m$
matrices. As usual, for $n\ge1$, let us denote by $\setA^{n}$
the set of all $n$-products of matrices from $\setA$;
$\setA^{0}=I$. For each $n\ge1$, define the quantity
\[
\bar{\rho}_{n}(\setA):=\max_{A_{i_{j}}\in\setA}
\rho(A_{i_{n}}\cdots
A_{i_{2}}A_{i_{1}}),
\]
where maximum is taken over all possible products of $n$
matrices from the set $\setA$, and $\rho(\cdot)$ denotes the
spectral radius of a matrix, that is the maximal magnitude of
its eigenvalues. Clearly, if $n > r$ then some matrices in the
product $A_{i_{n}}\cdots A_{i_{2}}A_{i_{1}}$ will occur several
times. The limit
\[
\bar{\rho}({\setA}):=
\limsup_{n\to\infty}\left(\bar{\rho}_{n}({\setA})\right)^{1/n}
\]
is called \emph{the generalized spectral radius} of the matrix
set $\setA$ \cite{DaubLag:LAA92,DaubLag:LAA01}.

Similarly, given a norm $\|\cdot\|$ in ${\mathbb{R}}^{m}$, for
each $n\ge1$, define the quantity
\[
\hat{\rho}_{n}({\setA}):= \max_{A_{i_{j}}\in\setA}
\|A_{i_{n}}\cdots
A_{i_{2}}A_{i_{1}}\|,
\]
where $\|A\|$, for a matrix $A$, is the matrix norm generated
by the vector norm $\|\cdot\|$ in ${\mathbb{R}}^{m}$, that is
$\|A\|=\sup_{\|x\|=1}\|Ax\|$. Then the limit
\begin{equation}\label{E-GenGF}
\hat{\rho}({\setA}):=
\limsup_{n\to\infty}\left(\hat{\rho}_{n}({\setA})\right)^{1/n}
\end{equation}
does not depend on the choice of the norm $\|\cdot\|$ and is
called \emph{the joint spectral radius}  of the matrix set
$\setA$ \cite{RotaStr:IM60}. When $r=1$ this definition
coincides with the famous Gelfand formula for the spectral
radius of a matrix \cite{Gelf:MatSb41:e} as in this case
$\setA=\{A\}$ is a singleton matrix set and
$\hat{\rho}_{n}({\setA})=\|A^{n}\|$.

For matrix sets $\setA$ consisting of a finite amount of
matrices, as is our case, the quantities $\bar{\rho}({\setA})$
and $\hat{\rho}({\setA})$ coincide with each other
\cite{BerWang:LAA92} and their common value is denoted as
\[
\rho({\setA}):=\bar{\rho}({\setA})=\hat{\rho}({\setA}),
\]
while the quantities $\bar{\rho}_{n}({\setA})$ and
$\hat{\rho}_{n}({\setA})$ form lower and upper bounds,
respectively, for the joint/generalized spectral radius:
\[
\bar{\rho}_{n}({\setA})\le
\bar{\rho}({\setA})=\hat{\rho}({\setA})\le
\hat{\rho}_{n}({\setA}),\qquad\forall~ n\ge0.
\]
This last formula may serve as a basis for a posteriori
estimating the accuracy of computation of $\rho({\setA})$. The
first algorithms of a kind in the context of control theory
problems have been suggested in \cite{BrayTong:TCS80}, for
linear inclusions in \cite{Bar:AIT88-2:e}, and for problems of
wavelet theory in
\cite{DaubLag:SIAMMAN92,DaubLag:LAA92,ColHeil:IEEETIT92}. Later
the computational efficiency of these algorithms was
essentially improved in \cite{Grip:LAA96,Maesumi:LAA96}.
Unfortunately, the common feature of all such algorithms is
that they do not provide any bounds for the number of
computational steps required to get desired accuracy of
approximation of $\rho({\setA})$.

Recently, in \cite{Koz:LAA09} explicit computable estimates for
the rate of convergence of the quantities $\|A^{n}\|^{1/n}$ to
$\rho(A)$ end their extension to the case of the joint spectral
radius were obtained. Probably, these results will help to make
more constructive the problem of evaluating of $\rho({\setA})$
by the generalized Gelfand formula \eqref{E-GenGF} .

Some works suggest different formulas to compute
$\rho({\setA})$. So, in \cite{ChenZhou:LAA00} it is shown that
\[
\rho({\setA})=
\limsup_{n\to\infty}\max_{A_{i_{j}}\in\setA}
\left|\tr(A_{i_{n}}\cdots
A_{i_{2}}A_{i_{1}})\right|^{1/n},
\]
where, as usual, $\tr(\cdot)$ denotes the trace of a matrix.

Given a norm $\|\cdot\|$ in ${\mathbb{R}}^{m}$, denote
\[
\|\setA\|:=\max_{A\in\setA}\|A\|.
\]
Then the spectral radius of the matrix set $\setA$ can be
defined by the equality
\begin{equation}\label{E-inf}
\rho({\setA})=\inf_{\|\cdot\|}\|\setA\|,
\end{equation}
where infimum is taken over all norms in ${\mathbb{R}}^{d}$
\cite{RotaStr:IM60,Els:LAA95}, and therefore
\[
\rho({\setA})\le\|\setA\|
\]
for any norm $\|\cdot\|$ in ${\mathbb{R}}^{m}$. For irreducible
matrix sets,\footnote{A matrix set $\setA$ is called
\emph{irreducible} if the matrices from $\setA$ have no common
invariant subspaces except $\{0\}$ and ${\mathbb{R}}^{m}$. In
\cite{KozPok:DAN92:e,KozPok:CADSEM96-005,KozPok:TRANS97} such a
matrix set was called quasi-controllable.} the infimum in
(\ref{E-inf}) is attained, and for such matrix sets there are
norms $\|\cdot\|$ in ${\mathbb{R}}^{m}$, called \emph{extremal
norms}, for which
\begin{equation}\label{E-extnorm}
    \|\setA\|\le\rho({\setA}).
\end{equation}

In analysis of the joint spectral radius ideas suggested by
N.E.~Barabanov \cite{Bar:AIT88-2:e,Bar:AIT88-3:e,Bar:AIT88-5:e}
play an important role. These ideas have got further
development in a variety of publications among which we would
like to distinguish \cite{Wirth:LAA02}.

\begin{theorem}[N.E.~Barabanov]
Let the matrix set $\setA=\{A_{1},\ldots,A_{r}\}$ be
irreducible. Then the
quantity $\rho$ is the joint (generalized) spectral radius of
the set $\setA$ iff there is a norm $\|\cdot\|$ in
${\mathbb{R}}^{m}$ such that
\begin{equation}\label{Eq-mane-bar}
\rho\|x\|\equiv
\max_{i}\|A_{i}x\|.
\end{equation}
\end{theorem}

Throughout the paper, a norm satisfying (\ref{Eq-mane-bar})
will be called a \emph{Barabanov norm} corresponding to the
matrix set $\setA$. Note that Barabanov norms are not unique.

Similarly, as is shown in \cite[Thm 3.3]{Prot:FPM96:e} and
\cite{Prot:FU98}, the value of $\rho$ equals to $\rho({\setA})$
if and only if for some central-symmetric convex
body\footnote{The set is called body if it contains at least
one interior point.} $S$ the following equality holds
\begin{equation}\label{E-protset}
    \rho S =\conv\left(\bigcup_{i=1}^{r}A_{i}S\right),
\end{equation}
where $\conv(\cdot)$ denotes the convex hull of a set and $\rho
S:=\{\rho x:\, x\in S\}$. As is noted in \cite{Prot:FPM96:e},
the relation (\ref{E-protset}) was proved by A.N. Dranishnikov
and S.V. Konyagin, so it is natural to call the
central-symmetric set $S$ the
\emph{Dranishnikov-Konyagin-Protasov set}. The set $S$ can be
treated as the unit ball of some norm $\|\cdot\|$ in
${\mathbb{R}}^{d}$ (recently this norm is usually called the
\emph{Protasov norm}). Note that Barabanov and Protasov norms
are the extremal norms, that is they satisfy the inequality
(\ref{E-extnorm}). In \cite{PWB:CDC05,Wirth:CDC05,PW:LAA08} it
is shown that Barabanov and Protasov norms are dual to each
other.

Remark that formulas (\ref{E-extnorm}), (\ref{Eq-mane-bar}) and
(\ref{E-protset}) define the joint or generalized spectral
radius for a matrix set in an apparently computationally
nonconstructive manner. In spite of that, namely such formulas
underlie quite a number of theoretical constructions (see,
e.g.,
\cite{Koz:CDC05:e,Koz:INFOPROC06:e,Wirth:CDC05,Wirth:LAA02,ParJdb:LAA08,Bar:CDC05})
and algorithms \cite{Prot:CDC05-1} for computation of
$\rho({\setA})$.

Different approaches for constructing Barabanov norms to
analyze properties of the joint (generalized) spectral radius
are discussed, e.g., in \cite{GugZen:LAA01,GugZen:LAA08} and
\cite[Sect. 6.6]{Theys:PhD05}. In
\cite{Koz:CDC05:e,Koz:INFOPROC06:e} the method of Barabanov
norms was the key instrument in disproving the Lagarias-Wang
Finiteness Conjecture. The related constructions were
essentially based on the study of the geometrical properties of
the unit balls of some specific Barabanov norms. In
\cite{GugZen:LAA08,GugZen:CDC05} the method of extremal
polytope norms was the key tool in investigation of the
finiteness properties of pairs of $2\times 2$ matrices.

In this context the situation when one fails to find among
current publications any detailed analysis of the geometrical
properties of the unit balls of Barabanov norms looks a bit
paradoxical. Partially this is explained by the fact that
Barabanov norms are defined nonconstructively, by an implicit
procedure. So, even in simplest cases it is very difficult to
visualize the shape of their unit balls. The present work may
be treated as the first step to make up this deficiency.

In the paper, an iteration procedure is considered that allows
to build numerically Barabanov norms for the irreducible matrix
sets and simultaneously to compute the joint spectral radius of
these sets. A similar iteration procedure is also discussed in
\cite{Koz:DAN09:e}.

The structure of the paper is as follows. In Introduction we
give basic definitions and present the motivation of the work.
In Section~\ref{S-iterscheme} the main iteration procedure is
introduced and Main Theorem stating convergence of this
procedure is formulated. The iteration procedure under
consideration is called the max-relaxation procedure since in
it the next approximation to the Barabanov norm is constructed
as the maximum of the current approximation and some auxiliary
norm. In Section~\ref{MTproof} proof Main Theorem is given. In
Section~\ref{S-comp2D}, to build simplest examples, the
max-relaxation scheme is adapted for computations with $2\times
2$ matrices. Results of numerical tests are illustrated by two
examples. At last, in concluding Section~\ref{S-rem} we discuss
some shortcomings of the proposed approach and formulate
further unresolved problems.

\section{Max-relaxation iteration scheme}\label{S-iterscheme}

Given $r,m\ge 1$, let $\setA=\{A_{1},\ldots,A_{r}\}$ be an
irreducible set of $m\times m$ real matrices.

Throughout the paper, a continuous function $\gamma(t,s)$
defined for $t,s> 0$ and satisfying
\[
\gamma(t,t)=t,\qquad
\min\{t,s\}<\gamma(t,s)<\max\{t,s\}\quad\textrm{for}~t\neq s,
\]
will be called \emph{an averaging
function}. Examples for averaging functions are:
\[
\gamma(t,s)=\frac{t+s}{2},\quad
\gamma(t,s)=\sqrt{ts},\quad
\gamma(t,s)=\frac{2ts}{t+s}.
\]

Given some averaging function $\gamma(\cdot,\cdot)$, a norm
$\|\cdot\|_{0}$ in ${\mathbb{R}}^{m}$, and  a vector
$e\in{\mathbb{R}}^{m}$ such that $\|e\|_{0}=1$, construct
recursively the norms $\|\cdot\|_{n}$ and
$\|\cdot\|^{\circ}_{n}$, $n=1,2,\ldots$, in accordance with the
following rules:
\medskip

\textbf{MR1:} \emph{if the norm $\|\cdot\|_{n}$ has been
already defined compute the quantities}
\begin{equation}\label{E-lohibounds}
 \rho^{+}_{n}=\max_{x\neq0}\frac{\max_{i}\|A_{i}x\|_{n}}{\|x\|_{n}},\quad
 \rho^{-}_{n}=\min_{x\neq0}\frac{\max_{i}\|A_{i}x\|_{n}}{\|x\|_{n}},\quad
 \gamma_{n}=\gamma(\rho^{-}_{n},\rho^{+}_{n});
\end{equation}

\textbf{MR2:} \emph{define the norms $\|\cdot\|_{n+1}$ and}
$\|\cdot\|^{\circ}_{n+1}$:
\begin{align}\label{E-auxnorm}
 \|x\|_{n+1}&=
 \max\left\{\|x\|_{n},
 ~\gamma^{-1}_{n}\max_{i}\|A_{i}x\|_{n}\right\},\\
\label{E-newnorm}
 \|x\|^{\circ}_{n+1}&=\|x\|_{n+1}/\|e\|_{n+1}.
\end{align}

Remark that the number of operations needed to perform one step
of algorithm \textbf{MR1-MR2} is of order
$rm^{2}\nu(\varepsilon)$, where $\nu(\varepsilon)$ is the
number of operations needed to compute, for an arbitrary vector
$x\in\mathbb{R}^{m}$, the value of the norm $\|x\|$ with a
relative accuracy $\varepsilon$. In general, the value
$\nu(\varepsilon)$ is of order $\varepsilon^{-m}$. So, the
total number of operations needed to perform $n$ steps of
algorithm \textbf{MR1-MR2} has the same rate of growth as
$nrm^{2}\varepsilon^{-m}$.

Remark also, that the procedure (\ref{E-lohibounds}) of
calculation of $\rho^{\pm}_{n}$ resembles the technique of
iterative approximation of the joint spectral radius suggested
in \cite{Grip:LAA96}.

\begin{main}
For any irreducible matrix set $\setA$, nonzero vector
$e\in{\mathbb{R}}^{m}$, initial norm $\|\cdot\|_{0}$, and any
averaging function $\gamma(t,s)$, the sequences
$\{\rho^{\pm}_{n}\}$ constructed by the iteration procedure
\textbf{MR1}, \textbf{MR2} converge to $\rho(\setA)$, and the
sequence of norms $\|\cdot\|_{n}$ converges uniformly on each
bounded set to some Barabanov norm $\|\cdot\|^{*}$ of the
matrix set $\setA$. Moreover, the sequence $\{\rho^{-}_{n}\}$
is nondecreasing, the sequence $\{\rho^{+}_{n}\}$ is
nonincreasing, and
$$
\rho^{-}_{n}\le \rho(\setA)\le
\rho^{+}_{n}
$$
for all $n=1,2,\ldots~$, which provides an a posteriori
estimate for the computational error of $\rho(\setA)$.
\end{main}

\section{Proof of Main Theorem}\label{MTproof}
Let us suppose that we managed to prove the following
assertions:
\begin{description}
\item[A1] \emph{the sequences $\{\rho^{+}_{n}\}$ and
    $\{\rho^{-}_{n}\}$ are convergent;}

\item[A2] \emph{the limits of the sequences
    $\{\rho^{+}_{n}\}$ and $\{\rho^{-}_{n}\}$ coincide:}
      \[
        \rho=\lim_{n\to\infty}\rho^{+}_{n}=\lim_{n\to\infty}\rho^{-}_{n};
      \]

\item[A3] \emph{the norms $\|\cdot\|^{\circ}_{n}$ converge
    pointwise to a limit $\|\cdot\|^{*}$.}
\end{description}

Then the function $\|\cdot\|^{*}$ will be a semi-norm in
${\mathbb{R}}^{m}$. Moreover, by (\ref{E-newnorm}) each norm
$\|\cdot\|^{\circ}_{n}$ meets the normalization condition
$\|e\|^{\circ}_{n}=1$, and hence
\[
 \|e\|^{*}=\lim_{n\to\infty}\|e\|^{\circ}_{n}= 1,
\]
which implies $\|x\|^{*}\not\equiv0$. Note also that due to
(\ref{E-newnorm}) the norms $\|\cdot\|^{\circ}_{n}$ differ from
$\|\cdot\|_{n}$ only by  numerical factors. Therefore, the
quantities  $\rho^{\pm}_{n}$ can be defined as
\begin{equation}\label{E-lohibounds-x}
 \rho^{+}_{n}=\max_{x\neq0}\frac{\max_{i}\|A_{i}x\|^{\circ}_{n}}{\|x\|^{\circ}_{n}},\quad
 \rho^{-}_{n}=\min_{x\neq0}\frac{\max_{i}\|A_{i}x\|^{\circ}_{n}}{\|x\|^{\circ}_{n}}.
\end{equation}
Then, passing to the limit in (\ref{E-lohibounds-x}), one can
conclude that the semi-norm $\|x\|^{*}$ satisfies the
\emph{Barabanov condition}
\[
\rho\|x\|^{*}=
\max_{i}\|A_{i}x\|^{*}.
\]
But as shown in \cite[Thm.~3]{Koz:INFOPROC06:e}, any semi-norm
$\|x\|^{*}\not\equiv0$ satisfying the Barabanov condition for
an irreducible matrix set is a Barabanov norm.

Thus, under assumptions \textbf{A1}, \textbf{A2} and
\textbf{A3}, the iteration procedure
(\ref{E-lohibounds})--(\ref{E-newnorm}) allows to build a
Barabanov norm and to find the joint spectral radius of the
matrix set $\setA$.

So, to complete the proof of Main Theorem we need to justify
assertions \textbf{A1}, \textbf{A2} and \textbf{A3} which will
be done in Sections~\ref{S-convnorm}--\ref{S-finA2}. In
Section~\ref{S-convnorm} we establish convergence of the
sequence of norms $\{\|\cdot\|^{\circ}_{n}\}$ to some limit
which allows to prove in Lemma~\ref{L-A3true} that Assertion
\textbf{A3} is a corollary of Assertions \textbf{A1} and
\textbf{A2}. Section~\ref{S-rhorel} demonstrates that the
quantities $\{\rho^{-}_{n}\}$ form the family of lower bounds
for the joint spectral radius $\rho$ of the matrix set $\setA$,
while the quantities $\{\rho^{+}_{n}\}$ form the family of
upper bounds for $\rho$.  In Section~\ref{S-convrho} we prove
that the sequences $\{\rho^{\pm}_{n}\}$ are bounded and
monotone which implies the existence of the limits
$\rho^{-}=\lim_{n\to\infty}\rho^{-}_{n}$ and
$\rho^{+}=\lim_{n\to\infty}\rho^{+}_{n}$. At last, in
Sections~\ref{S-newnorms} and \ref{S-omega} we prove that
$\rho^{-}=\rho^{+}$ which allows to justify in
Section~\ref{S-finA2} the validity of Assertions \textbf{A1}
and \textbf{A2} and thus to finalize the proof of Main Theorem.

\subsection{Convergence of the sequence of norms
$\{\|\cdot\|^{\circ}_{n}\}$}\label{S-convnorm}

Given a pair of norms $\|\cdot\|'$ and $\|\cdot\|''$ in
${\mathbb{R}}^{m}$ define the quantities
\begin{equation}\label{E-eccentr}
    e^{-}(\|\cdot\|',\|\cdot\|'')=\min_{x\neq0}\frac{\|x\|'}{\|x\|''},\quad
    e^{+}(\|\cdot\|',\|\cdot\|'')=\max_{x\neq0}\frac{\|x\|'}{\|x\|''}.
\end{equation}

Since all norms in ${\mathbb{R}}^{m}$ are equivalent to each
other, the quantities $e^{-}(\|\cdot\|',\|\cdot\|'')$ and
$e^{+}(\|\cdot\|',\|\cdot\|'')$ are correctly defined and
\[
0< e^{-}(\|\cdot\|',\|\cdot\|'')\le
e^{+}(\|\cdot\|',\|\cdot\|'')< \infty.
\]
Therefore the quantity
\begin{equation}\label{E-defeccentr}
\ecc(\|\cdot\|',\|\cdot\|'')=
\frac{e^{+}(\|\cdot\|',\|\cdot\|'')}{e^{-}(\|\cdot\|',\|\cdot\|'')}\ge 1,
\end{equation}
which is called \textit{the eccentricity} of the norm $\|\cdot\|'$
with respect to the norm $\|\cdot\|''$ (see, e.g., \cite{Wirth:CDC05}), is also correctly defined.

\begin{lemma}\label{L-eccbound}
Let $\|\cdot\|^{*}$ be a Barabanov norm for the matrix set
$\setA$. Then
\begin{equation}\label{E-eqecc}
\ecc(\|\cdot\|^{\circ}_{n},\|\cdot\|^{*})=
\ecc(\|\cdot\|_{n},\|\cdot\|^{*}),\quad\forall~n,
\end{equation}
and the sequence of the numbers $\ecc(\|\cdot\|_{n},\|\cdot\|^{*})$ is nonincreasing.
\end{lemma}

\proof Note first that by (\ref{E-newnorm}) each norm $\|\cdot\|^{\circ}_{n}$ differs from the corresponding norm  $\|\cdot\|_{n}$ only by a numerical factor. From this, by the definition (\ref{E-eccentr}), (\ref{E-defeccentr}) of the eccentricity of one norm with respect to another, the equality (\ref{E-eqecc}) follows.

Denote by $\rho$ the joint spectral radius of the matrix set $\setA$. Then, by definitions of the function $e^{+}(\cdot)$ and of the Barabanov norm $\|\cdot\|^{*}$, from
(\ref{E-lohibounds}), (\ref{E-auxnorm}) we obtain:
\begin{multline*}
\|x\|_{n+1}=\max\left\{\|x\|_{n},
 ~\gamma^{-1}_{n}\max_{i}\|A_{i}x\|_{n}\right\}\le\\
 \le e^{+}(\|\cdot\|_{n},\|\cdot\|^{*})  \max\left\{\|x\|^{*},
 ~\gamma^{-1}_{n}\max_{i}\|A_{i}x\|^{*}\right\}=\\
= e^{+}(\|\cdot\|_{n},\|\cdot\|^{*}) \max\left\{\|x\|^{*},
 ~\gamma^{-1}_{n}\rho\|x\|^{*}\right\}.
\end{multline*}
Therefore
\begin{equation}\label{E-eccupbound}
   e^{+}(\|\cdot\|_{n+1},\|\cdot\|^{*})\le
   e^{+}(\|\cdot\|_{n},\|\cdot\|^{*})
\max\left\{1,
 ~\gamma^{-1}_{n}\rho\right\}.
\end{equation}

Similarly, by definitions of the function $e^{-}(\cdot)$ and of the Barabanov norm  $\|\cdot\|^{*}$, from (\ref{E-lohibounds}),
(\ref{E-auxnorm}) we obtain:
\begin{multline*}
\|x\|_{n+1}=  \max\left\{\|x\|_{n},
 ~\gamma^{-1}_{n}\max_{i}\|A_{i}x\|_{n}\right\}\ge\\
 \ge e^{-}(\|\cdot\|_{n},\|\cdot\|^{*})  \max\left\{\|x\|^{*},
 ~\gamma^{-1}_{n}\max_{i}\|A_{i}x\|^{*}\right\}=\\
= e^{-}(\|\cdot\|_{n},\|\cdot\|^{*}) \max\left\{\|x\|^{*},
 ~\gamma^{-1}_{n}\rho\|x\|^{*}\right\}.
\end{multline*}
Therefore
\begin{equation}\label{E-ecclobound}
   e^{-}(\|\cdot\|_{n+1},\|\cdot\|^{*})\ge
   e^{-}(\|\cdot\|_{n},\|\cdot\|^{*})
\max\left\{1,
 ~\gamma^{-1}_{n}\rho\right\}.
\end{equation}

By dividing termwise the inequality (\ref{E-eccupbound}) on
(\ref{E-ecclobound}) we get
\[
\ecc(\|\cdot\|_{n+1},\|\cdot\|^{*})=
\frac{e^{+}(\|\cdot\|_{n+1},\|\cdot\|^{*})}{e^{-}(\|\cdot\|_{n+1},\|\cdot\|^{*})}\le
\frac{e^{+}(\|\cdot\|_{n},\|\cdot\|^{*})}{e^{-}(\|\cdot\|_{n},\|\cdot\|^{*})} =
\ecc(\|\cdot\|_{n},\|\cdot\|^{*}).
\]
Hence, the sequence $\{\ecc(\|\cdot\|_{n},\|\cdot\|^{*})\}$ is nonincreasing. \qed

Denote by $N_{\mathrm{loc}}(\mathbb{R}^{m})$ the topological
space of norms in $\mathbb{R}^{m}$ with the topology of uniform
convergence on bounded subsets of $\mathbb{R}^{m}$.

\begin{corollary}\label{C1-L-eccbound}
The sequence of norms $\{\|\cdot\|^{\circ}_{n}\}$ is compact in
$N_{\mathrm{loc}}(\mathbb{R}^{m})$.
\end{corollary}

\proof For each $n$ and any $x\neq0$, by the definition
(\ref{E-eccentr}) of the functions $e^{+}(\cdot)$ and
$e^{-}(\cdot)$ the following relations hold
\[
e^{-}(\|\cdot\|_{n},\|\cdot\|^{*})\le\frac{\|x\|_{n}}{\|x\|^{*}}\le e^{+}(\|\cdot\|_{n},\|\cdot\|^{*}),
\]
and then
\[
e^{-}(\|\cdot\|_{n},\|\cdot\|^{*})\le\frac{\|e\|_{n}}{\|e\|^{*}}\le e^{+}(\|\cdot\|_{n},\|\cdot\|^{*}),
\]
from which
\[
\frac{1}{\ecc(\|\cdot\|_{n},\|\cdot\|^{*})}\frac{\|x\|^{*}}{\|e\|^{*}}
\|e\|_{n}\le\|x\|_{n}\le \ecc(\|\cdot\|_{n},\|\cdot\|^{*})
\frac{\|x\|^{*}}{\|e\|^{*}}\|e\|_{n}.
\]
Since here by construction the norms $\{\|\cdot\|^{\circ}_{n}\}$ satisfy the normalization condition $\|e\|^{\circ}_{n}\equiv 1$, and  by Lemma~\ref{L-eccbound} $\ecc(\|\cdot\|^{\circ}_{n},\|\cdot\|^{*})\le
\ecc(\|\cdot\|^{\circ}_{0},\|\cdot\|^{*})$, we have
\[
\frac{1}{\ecc(\|\cdot\|_{0},\|\cdot\|^{*})}\frac{\|x\|^{*}}{\|e\|^{*}}\le\|x\|_{n}\le
\ecc(\|\cdot\|_{0},\|\cdot\|^{*})\frac{\|x\|^{*}}{\|e\|^{*}}.
\]

Therefore the norms $\|\cdot\|_{n}$, $n\ge 1$, are
equicontinuous and uniformly bounded on each bounded subset of
$\mathbb{R}^{m}$. Moreover, their values are also uniformly
separated from zero on each bounded subset of $\mathbb{R}^{m}$
separated from zero. From here by the Arzela-Ascoli theorem the
statement of the corollary follows. \qed

\begin{corollary}\label{C2-L-eccbound}
If at least one of subsequences of norms from
$\{\|\cdot\|^{\circ}_{n}\}$ converges in
$N_{\mathrm{loc}}(\mathbb{R}^{m})$ to some Barabanov norm then
the whole sequence $\{\|\cdot\|^{\circ}_{n}\}$ also converges
in $N_{\mathrm{loc}}(\mathbb{R}^{m})$ to the same Barabanov
norm.
\end{corollary}

\proof Let $\{\|\cdot\|^{\circ}_{n_{k}}\}$ be a subsequence of
$\{\|\cdot\|^{\circ}_{n}\}$ which converges in
$N_{\mathrm{loc}}(\mathbb{R}^{m})$ to some Barabanov norm
$\|\cdot\|^{*}$. Then by definition of the eccentricity of one
norm with respect to another
\[
\ecc(\|\cdot\|^{\circ}_{n_{k}},\|\cdot\|^{*})\to 1\quad\textrm{as~}k\to\infty.
\]
Here by Lemma~\ref{L-eccbound} the eccentricities
$\ecc(\|\cdot\|^{\circ}_{n},\|\cdot\|^{*})$ are nonincreasing in $n$, and then the following stronger relation holds
\begin{equation}\label{E-eccconv}
    \ecc(\|\cdot\|^{\circ}_{n},\|\cdot\|^{*})\to 1\quad\textrm{as~}n\to\infty.
\end{equation}

Note now that by the definition (\ref{E-eccentr}),
(\ref{E-defeccentr}) of the eccentricity of one norm with respect to another
\[
\frac{1}{\ecc(\|\cdot\|^{\circ}_{n},\|\cdot\|^{*})} \le
\frac{\|x\|^{\circ}_{n}}{\|x\|^{*}}\le \ecc(\|\cdot\|^{\circ}_{n},\|\cdot\|^{*}),
\]
from which by (\ref{E-eccconv}) it follows that the sequence of norms $\{\|\cdot\|^{\circ}_{n}\}$ converges in  space $N_{\mathrm{loc}}(\mathbb{R}^{m})$ to the norm $\|\cdot\|^{*}$. \qed

\begin{lemma}\label{L-A3true}
Assertion \textbf{A3} is a corollary of Assertions \textbf{A1} and \textbf{A2}.
\end{lemma}

\proof By Corollary~\ref{C1-L-eccbound} the sequence of norms
$\{\|\cdot\|^{\circ}_{n}\}$ has a subsequence
$\{\|\cdot\|^{\circ}_{n_{k}}\}$ that converges in space $N_{\mathrm{loc}}(\mathbb{R}^{m})$ to some norm
$\|\cdot\|^{*}$. Then, passing to the limit in (\ref{E-lohibounds-x}) as $n=n_{k}\to\infty$, we get by Assertions \textbf{A1} and \textbf{A2}:
\[
\rho=\frac{\max_{i}\|A_{i}x\|^{*}}{\|x\|^{*}},\quad \forall~x\neq0,
\]
which means that $\|\cdot\|^{*}$ is a Barabanov norm for the
matrix set $\setA$. This and Corollary~\ref{C2-L-eccbound} then
imply that the sequence $\{\|\cdot\|^{\circ}_{n}\}$ converges
in space $N_{\mathrm{loc}}(\mathbb{R}^{m})$ to the Barabanov
norm $\|\cdot\|^{*}$. Assertion \textbf{A3} is proved. \qed

So, in view of Lemma~\ref{L-A3true} to prove that the iteration
procedure (\ref{E-lohibounds})--(\ref{E-newnorm}) is convergent
it suffices to verify only that Assertions \textbf{A1} and
\textbf{A2} hold.

\subsection{Relations between
$\rho^{\pm}_{n}$ and $\rho$}\label{S-rhorel}

The following lemma provides a way to estimate the spectral
radius of a matrix set.

\begin{lemma}\label{L-rhorel}
Let $\alpha,\beta$ be numbers such that in some norm
$\|\cdot\|$ the inequalities
\[
\alpha\|x\|\le\max_{A_{i}\in\setA}\|A_{i}x\|\le\beta\|x\|,
\]
hold. Then $\alpha\le\rho\le\beta$, where $\rho$ is the joint
spectral radius of the matrix set $\setA$.
\end{lemma}

\proof Let $\|\cdot\|^{*}$ be some Barabanov norm for the
matrix set $\setA$. Since all norms in ${\mathbb{R}}^{m}$ are
equivalent, there are constants $\sigma^{-}>0$ and
$\sigma^{+}<\infty$ such that
\begin{equation}\label{E-sigmaest}
    \sigma^{-}\|x\|^{*}\le\|x\|\le\sigma^{+}\|x\|^{*}.
\end{equation}
Consider, for each $k=1,2,\ldots$, the functions
\[
\Delta_{k}(x)=
\max_{1\le i_{1},i_{2},\ldots,i_{k}\le r}\|A_{i_{k}}\dots A_{i_{2}}A_{i_{1}}x\|.
\]
Then, as is easy to see,
\begin{equation}\label{E-Dbounds}
\alpha^{k}\|x\|\le\Delta_{k}(x)\le\beta^{k}\|x\|.
\end{equation}

Similarly consider, for each $k=1,2,\ldots$, the functions
\[
\Delta^{*}_{k}(x)=
\max_{1\le i_{1},i_{2},\ldots,i_{k}\le r}\|A_{i_{k}}\dots A_{i_{2}}A_{i_{1}}x\|^{*}.
\]
For these functions, by definition of Barabanov norms the
following identity hold
\begin{equation}\label{E-Dstar}
\Delta^{*}_{k}(x)\equiv\rho^{k}\|x\|^{*},
\end{equation}
which is stronger than (\ref{E-Dbounds}).

Now, note that (\ref{E-sigmaest}) and the definition of the
functions $\Delta_{k}(x)$ and $\Delta^{*}_{k}(x)$ imply
\[
\sigma^{-}\Delta^{*}_{k}(x)\le\Delta_{k}(x)\le\sigma^{+}\Delta^{*}_{k}(x).
\]
Then, by (\ref{E-Dbounds}), (\ref{E-Dstar}),
\[
\frac{\sigma^{-}}{\sigma^{+}}\alpha^{k}
\le\rho^{k}\le\frac{\sigma^{+}}{\sigma^{-}}\beta^{k},\quad
\forall~k,
\]
from which the required estimates $\alpha\le\rho\le\beta$
follow. \qed

So, Lemma~\ref{L-rhorel} and the definition
(\ref{E-lohibounds}) of $\rho^{\pm}_{n}$ imply that the
quantities $\{\rho^{-}_{n}\}$ form the family of lower bounds
for the joint spectral radius $\rho$ of the matrix set $\setA$,
while the quantities $\{\rho^{+}_{n}\}$ form the family of
upper bounds for $\rho$. This allows to estimate a posteriori
errors of computation of the joint spectral radius with the
help of the iteration procedure
(\ref{E-lohibounds})--(\ref{E-newnorm}).

\subsection{Convergence of the sequences
$\{\rho^{\pm}_{n}\}$}\label{S-convrho}

Estimate the value of $\max_{i}\|A_{i}x\|_{n+1}$. By definition,
\begin{multline*}
\max_{i}\|A_{i}x\|_{n+1}
=\max_{i}\left\{\max\left\{\|A_{i}x\|_{n},
 ~\gamma^{-1}_{n}\max_{j}\|A_{i}A_{j}x\|_{n}\right\}\right\}=\\
= \max\left\{\max_{i}\|A_{i}x\|_{n},
~\gamma^{-1}_{n}\max_{j}\max_{i}\|A_{i}A_{j}x\|_{n}\right\}.
\end{multline*}
Here by the definition (\ref{E-lohibounds}) of the quantities
$\rho^{\pm}_{n}$ the right-hand part of the chain of equalities can be estimated as follows:
\begin{multline*}
\rho^{-}_{n}\max\left\{\|x\|_{n},
~\gamma^{-1}_{n}\max_{j}\|A_{j}x\|_{n}\right\}\le\\
\le \max\left\{\max_{i}\|A_{i}x\|_{n},
~\gamma^{-1}_{n}\max_{j}\max_{i}\|A_{i}A_{j}x\|_{n}\right\}\le\\
\rho^{+}_{n}\max\left\{\|x\|_{n},
~\gamma^{-1}_{n}\max_{j}\|A_{j}x\|_{n}\right\}.
\end{multline*}
Therefore, by definition of the norm $\|x\|_{n+1}$,
\[
\rho^{-}_{n}\|x\|_{n+1}\le\max_{i}\|A_{i}x\|_{n+1}\le
\rho^{+}_{n}\|x\|_{n+1},
\]
from which
\[
\rho^{-}_{n}\le\frac{\max_{i}\|A_{i}x\|_{n+1}}{\|x\|_{n+1}}\le
\rho^{+}_{n},\quad\forall~x\neq0,
\]
and then,
\[
\rho^{-}_{n}\le\rho^{-}_{n+1}\le
\rho^{+}_{n+1}\le
\rho^{+}_{n}.
\]
So, the following lemma holds.
\begin{lemma}\label{L-monrho}
The sequence $\{\rho^{-}_{n}\}$ is bounded from above by each member of the sequence $\{\rho^{+}_{n}\}$ and is nondecreasing. The sequence $\{\rho^{+}_{n}\}$ is bounded from below by each member of the sequence  $\{\rho^{-}_{n}\}$ and is nonincreasing.
\end{lemma}

In view of Lemma~\ref{L-monrho} there are the limits
\[
\rho^{-}=\lim_{n\to\infty}\rho^{-}_{n},\quad
\rho^{+}=\lim_{n\to\infty}\rho^{+}_{n},\quad
\gamma=\lim_{n\to\infty}\gamma_{n}=
\lim_{n\to\infty}\gamma(\rho^{-}_{n},\rho^{+}_{n}),
\]
where
\[
\rho^{-}\le\gamma\le\rho^{+},
\]
which means that Assertion \textbf{A1} holds. Hence, to prove
that the iteration procedure
(\ref{E-lohibounds})--(\ref{E-newnorm}) is convergent it
remains only to justify Assertion \textbf{A2}: $\rho^{-}=
\rho^{+}$.

To prove that $\rho^{-}= \rho^{+}$, below it will be supposed
the contrary, which will lead us to a contradiction.

\subsection{Transition to a new sequence of norms}\label{S-newnorms}

To simplify further constructions we will switch over to a new
sequence of norms for which the quantities $\rho^{\pm}_{n}$
will be independent of $n$.

By Corollary~\ref{C1-L-eccbound} the sequence of the norms
$\|\cdot\|^{\circ}_{n}$ is compact in space
$N_{\mathrm{loc}}(\mathbb{R}^{m})$. Hence, there is a
subsequence of indices $\{n_{k}\}$ such that the the norms
$\|\cdot\|^{\circ}_{n_{k}}=\|\cdot\|_{n_{k}}/\|e\|_{n_{k}}$
converge to some norm $\|\cdot\|^{\bullet}_{0}$ satisfying the
normalization condition $\|e\|^{\bullet}_{0}=1$. Then, passing
to the limit in (\ref{E-lohibounds-x}), by Lemma~\ref{L-monrho}
we obtain:
\[
\rho^{+}=\max_{x\neq0}\frac{\max_{i}\|A_{i}x\|^{\bullet}_{0}}{\|x\|^{\bullet}_{0}},\quad
\rho^{-}=\min_{x\neq0}\frac{\max_{i}\|A_{i}x\|^{\bullet}_{0}}{\|x\|^{\bullet}_{0}},\quad
\gamma=\gamma(\rho^{-},\rho^{+}).
\]
Now by induction the following statement can be easily proved.
\begin{lemma}\label{L-seqsharp}
For each $n=0,1,2,\ldots$, the sequence of the norms
$\|\cdot\|_{n_{k}+n}/\|e\|_{n_{k}}$ converges to some norm
$\|\cdot\|^{\bullet}_{n}$. Moreover, for each $n=0,1,2,\ldots$,
we have the equalities
\begin{equation}\label{E-lohisharp}
 \max_{x\neq0}\frac{\max_{i}\|A_{i}x\|^{\bullet}_{n}}{\|x\|^{\bullet}_{n}}= \rho^{+},\quad
 \min_{x\neq0}\frac{\max_{i}\|A_{i}x\|^{\bullet}_{n}}{\|x\|^{\bullet}_{n}}= \rho^{-},
\end{equation}
and the recurrent relations
\begin{equation}\label{E-recsharp}
 \|x\|^{\bullet}_{n+1}=
 \max\left\{\|x\|^{\bullet}_{n},
 ~\gamma^{-1} \max_{i}\|A_{i}x\|^{\bullet}_{n}\right\}.
\end{equation}
\end{lemma}

\subsection{Sets $\omega_{n}$}\label{S-omega}

Define, for each $n=0,1,2,\ldots$, the set
\begin{equation}\label{E-defomega}
    \omega_{n}=\left\{x\in\mathbb{R}^{m}:~\rho^{-}\|x\|^{\bullet}_{n}=
    \max_{i}\|A_{i}x\|^{\bullet}_{n}\right\}.
\end{equation}
By (\ref{E-lohisharp}) $\omega_{n}$ is the set on which the quantity
\[
\frac{\max_{i}\|A_{i}x\|^{\bullet}_{n}}{\|x\|^{\bullet}_{n}}
\]
attains its minimum.

\begin{lemma}\label{L-eqnonomega}
If $x\in\omega_{n}$ then
$\|x\|^{\bullet}_{n+1}=\|x\|^{\bullet}_{n}$.
\end{lemma}

\proof The statement of the lemma is obvious for $x=0$. So,
suppose that $x\in\omega_{n}$, $x\neq 0$. In this case
(\ref{E-defomega}) and the inequalities $\rho^{-}\le\rho^{+}$
imply
\[
\max_{i}\|A_{i}x\|^{\bullet}_{n}=
\rho^{-}\|x\|^{\bullet}_{n}\le \gamma \|x\|^{\bullet}_{n}
\]
or, what is the same,
\[
\|x\|^{\bullet}_{n}\ge\gamma^{-1}
 \max_{i}\|A_{i}x\|^{\bullet}_{n}.
\]
From here by the definition (\ref{E-recsharp}) of the norm
$\|\cdot\|^{\bullet}_{n+1}$ we get the required equality:
\[
\|x\|^{\bullet}_{n+1}=
 \max\left\{\|x\|^{\bullet}_{n},
 ~\gamma^{-1}
 \max_{i}\|A_{i}x\|^{\bullet}_{n}\right\}=
\|x\|^{\bullet}_{n}.
\]
The lemma is proved. \qed

\begin{lemma}\label{L-omegadec}
If $\rho^{-}<\rho^{+}$ then $\omega_{n+1}\subseteq\omega_{n}$
for each $n=0,1,2,\ldots$.
\end{lemma}

\proof Let $x\in\omega_{n+1}$. If $x=0$ then clearly
$x\in\omega_{n}$, so suppose that $x\neq 0$. By definitions of
the set $\omega_{n+1}$ and of the norm
$\|\cdot\|^{\bullet}_{n}$ the following equalities take place:
\begin{multline}\label{E-om}
\max_{i}\|A_{i}x\|^{\bullet}_{n+1}
=\max_{i}\left\{\max\left\{\|A_{i}x\|^{\bullet}_{n},
 ~\gamma^{-1} \max_{j}\|A_{j}A_{i}x\|^{\bullet}_{n}\right\}\right\}=\\
=\max\left\{\max_{i}\|A_{i}x\|^{\bullet}_{n},
 ~\gamma^{-1} \max_{i,j}\|A_{j}A_{i}x\|^{\bullet}_{n}\right\}=\\
=
\rho^{-}\|x\|^{\bullet}_{n+1}=\rho^{-}\max\left\{\|x\|^{\bullet}_{n},
~\gamma^{-1} \max_{i}\|A_{i}x\|^{\bullet}_{n}\right\}.
\end{multline}
Let here
\begin{equation}\label{E-wrongsuppose}
\|x\|^{\bullet}_{n}\le
\gamma^{-1} \max_{i}\|A_{i}x\|^{\bullet}_{n}.
\end{equation}
Then from (\ref{E-om}) it follows that
\[
\max\left\{\max_{i}\|A_{i}x\|^{\bullet}_{n},
 ~\gamma^{-1} \max_{i,j}\|A_{j}A_{i}x\|^{\bullet}_{n}\right\}
= \rho^{-}\|x\|^{\bullet}_{n+1}=
~\gamma^{-1}\rho^{-}\max_{i}\|A_{i}x\|^{\bullet}_{n}.
\]
But by the conditions of the lemma  $\rho^{-}<\rho^{+}$. Then
$\gamma=\gamma(\rho^{-},\rho^{+})>\rho^{-}$, and the right-hand part of the above equalities is strictly less than
$\max_{i}\|A_{i}x\|^{\bullet}_{n}$. A contradiction, since the left-hand part of the same equalities is no less than
$\max_{i}\|A_{i}x\|^{\bullet}_{n}$.

The above contradiction is caused by the assumption
(\ref{E-wrongsuppose}), and therefore it is proved that the
condition $x\neq0\in\omega_{n+1}$ implies the strict inequality
\[
\|x\|^{\bullet}_{n}>
\gamma^{-1} \max_{i}\|A_{i}x\|^{\bullet}_{n}.
\]
In this case from (\ref{E-om}) it follows that
\begin{equation}\label{E-maineq}
\max\left\{\max_{i}\|A_{i}x\|^{\bullet}_{n},
~\gamma^{-1} \max_{i,j}\|A_{j}A_{i}x\|^{\bullet}_{n}\right\}=
\rho^{-}\|x\|^{\bullet}_{n}.
\end{equation}

Let us show that the equality (\ref{E-maineq}) implies
\begin{equation}\label{E-finomega}
    \max_{i}\|A_{i}x\|^{\bullet}_{n}=\rho^{-}\|x\|^{\bullet}_{n}.
\end{equation}
Indeed, supposing the contrary, by definition of the quantity $\rho^{-}$, there should be valid the strict inequality $\max_{i}\|A_{i}x\|^{\bullet}_{n}>\rho^{-}\|x\|^{\bullet}_{n}$.
Then the left-hand part of the equality (\ref{E-maineq}) should be strictly greater than $\rho^{-}\|x\|^{\bullet}_{n}$, that is greater than the right-hand part of the same equality, which is impossible. This last contradiction shows that the equality
(\ref{E-finomega}) holds as soon as $x\neq0\in\omega_{n+1}$, which means by (\ref{E-lohisharp}) that $x\in\omega_{n}$. \qed

\begin{corollary}\label{C-mainomegaprop}
If $\rho^{-}<\rho^{+}$ then
$\omega:=\cap_{n\ge0}\omega_{n}\neq0$ and
\begin{equation}\label{E-eqsharpnorms}
\|x\|^{\bullet}_{0}=\|x\|^{\bullet}_{1}=\dots=\|x\|^{\bullet}_{n}=\dots,
\quad\forall~x\neq0\in\omega.
\end{equation}
\end{corollary}

\proof By Lemma~\ref{L-omegadec} $\{\omega_{n}\}$ is the family
of embedded closed conic\footnote{A set $X$ is called conic if
together with each its point $x$ it contains the ray $\{tx:~
t\ge0\}$.} sets. Then the intersection $\omega$ of these sets
is also a closed conic set such that $\omega\neq\{0\}$.

By definition of the set $\omega$, if $x\in\omega$ then
$x\in\omega_{n}$ for every integer $n\ge0$. Hence, by
Lemma~\ref{L-eqnonomega}
$\|x\|^{\bullet}_{n+1}=\|x\|^{\bullet}_{n}$, from which the
equalities (\ref{E-eqsharpnorms}) follow. \qed

\subsection{Completion of the proof of Assertion A2}\label{S-finA2}

By Corollary~\ref{C-mainomegaprop}
there is a non-zero vector $g$ on which all the norms
$\|\cdot\|^{\bullet}_{n}$ take the same values:
\[
\|g\|^{\bullet}_{0}
=\|g\|^{\bullet}_{1}=\dots=\|g\|^{\bullet}_{n}=\cdots.
\]
Then, due to uniform boundedness of the eccentricities of the
norms $\|\cdot\|^{\bullet}_{n}$ with respect to some Barabanov
norm $\|\cdot\|^{*}$ (this fact can be proved by verbatim
repetition of the analogous proof for the norms
$\|\cdot\|_{n}$), the norms $\|\cdot\|^{\bullet}_{n}$ form a
family which is uniformly bounded and equicontinuous with
respect to the Barabanov norm $\|\cdot\|^{*}$:
\[
\exists~ \mu^{\pm}\in(0,\infty):\quad
\mu^{-}\|x\|^{*}\le\|x\|^{\bullet}_{n}\le\mu^{+}\|x\|^{*},\qquad n=0,1,2,\ldots.
\]
Hence, by the Arzela-Ascoli theorem the family of norms $\{\|\cdot\|^{\bullet}_{n}\}$ is compact in $N_{\mathrm{loc}}(\mathbb{R}^{m})$.

From the definition (\ref{E-recsharp}) of the norms
$\|\cdot\|^{\bullet}_{n}$ it follows also that
\[
 \|x\|^{\bullet}_{n+1}=
 \max\left\{\|x\|^{\bullet}_{n},
 ~\gamma^{-1}
 \max_{i}\|A_{i}x\|^{\bullet}_{n}\right\}
 \ge\|x\|^{\bullet}_{n}.
\]
Then the norms $\|\cdot\|^{\bullet}_{n}$ are monotone
increasing in $n$ and bounded (with respect to the Barabanov
norm $\|\cdot\|^{*}$) and therefore they pointwise converge to
some norm $\|\cdot\|^{\bullet}$. Moreover, since the family of
norms $\{\|\cdot\|^{\bullet}_{n}\}$ is equicontinuous with
respect to the Barabanov norm $\|\cdot\|^{*}$, the norms
$\|\cdot\|^{\bullet}_{n}$ converge to the norm
$\|\cdot\|^{\bullet}$ in space
$N_{\mathrm{loc}}(\mathbb{R}^{m})$.

Now, passing to the limit in the relations
\[
 \|x\|^{\bullet}_{n+1}
= \max\left\{\|x\|^{\bullet}_{n},
 ~\gamma^{-1}
 \max_{i}\|A_{i}x\|^{\bullet}_{n}\right\}
 \ge\gamma^{-1}
 \max_{i}\|A_{i}x\|^{\bullet}_{n},
\]
which follow from (\ref{E-recsharp}), we obtain
\[
 \|x\|^{\bullet}
 \ge\gamma^{-1}
 \max_{i}\|A_{i}x\|^{\bullet}.
\]
From here
\begin{equation}\label{E-finestim}
\max_{x\neq0}
\frac{\max_{i}\|A_{i}x\|^{\bullet}}%
{\|x\|^{\bullet}}\le\gamma.
\end{equation}

On the other hand, passing to the limit in the first relation of
(\ref{E-lohisharp}), we obtain
\begin{equation}\label{E-finestim1}
\max_{x\neq0}
\frac{\max_{i}\|A_{i}x\|^{\bullet}}%
{\|x\|^{\bullet}}=\rho^{+}.
\end{equation}

Relations (\ref{E-finestim}) and (\ref{E-finestim1}) imply the inequality $\rho^{+}\le\gamma$ which contradicts the assumption
$\rho^{-}< \rho^{+}$ because by definition of the function $\gamma(\cdot,\cdot)$ the condition $\rho^{-}<
\rho^{+}$ implies the inequality
$\gamma=\gamma(\rho^{-},\rho^{+})<\rho^{+}$.

The obtained contradiction completes the proof of the equality $\rho^{-}= \rho^{+}$ as well as of the convergence of the iteration procedure (\ref{E-lohibounds})--(\ref{E-newnorm}).

\section{Computational scheme for $2\times 2$ matrices}\label{S-comp2D}

Let $\setA=\{A_{1},\ldots,A_{r}\}$ be a set of real
$2\times 2$ matrices
\[
A_{i}=\left(\begin{array}{cc}
a^{(i)}_{11}&a^{(i)}_{12}\\[2mm]a^{(i)}_{21}&a^{(i)}_{22}
\end{array}\right).
\]

Let $(r,\varphi)$ be the polar coordinates in $\mathbb{R}^{2}$.
Then, for a vector $x\in\mathbb{R}^{2}$ with Cartesian
coordinates $x=\{x_{1},x_{2}\}$, we have
\[
x=\{r\cos\varphi,r\sin\varphi\}
\]
and
\[
r=r(x)=\sqrt{x_{1}^{2}+x_{2}^{2}},\quad
\varphi=\varphi(x)=\arctan\left(x_{2}/x_{1}\right).
\]

\begin{figure}[!htbp]
\begin{center}
\includegraphics[scale=0.8,clip]{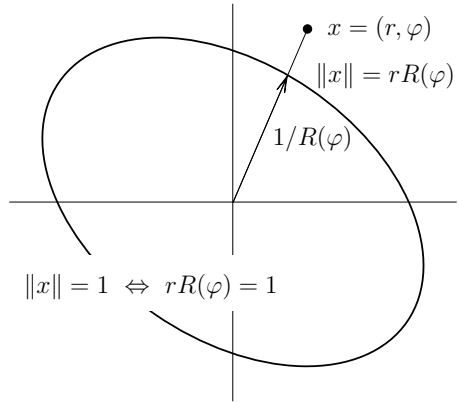}
\caption{Definition of the function $R(\varphi)$.}\label{F-polcoord}
\end{center}
\end{figure}

Define, for an arbitrary norm $\|\cdot\|$, the function
\[
R(\varphi)=\|\{\cos\varphi,\sin\varphi\}\|.
\]
Then the norm $\|x\|$ of the vector $x$ with polar coordinates
$(r,\varphi)$ is determined by the equality
\begin{equation}\label{E-rR}
\|x\|=rR(\varphi),
\end{equation}
and the unit sphere in the norm $\|\cdot\|$ is determined as
the geometrical locus of the vectors $x$ polar coordinates of
which satisfy the relations
\[
rR(\varphi)\equiv 1\quad\textrm{or}\quad r=\frac{1}{R(\varphi)},
\]
see Fig.~\ref{F-polcoord}.

Now, let $R_{n}(\varphi)$ be the function defining in the polar
coordinates the graph of the unit sphere $\|x\|_{n}=1$ of the
norm $\|\cdot\|_{n}$ determined by the iteration procedure
(\ref{E-lohibounds})--(\ref{E-newnorm}). Rewrite the relations
(\ref{E-lohibounds})--(\ref{E-newnorm}) in terms of the
functions $R_{n}(\varphi)$. To do it we should express the
quantities $\|A_{i}x\|_{n}$, $i=0,2,\ldots,r$, in terms of the
functions $R_{n}(\varphi)$.

By (\ref{E-rR})
\[
\|A_{i}x\|_{n}=r(A_{i}x)R_{n}(\varphi(A_{i}x)).
\]
Here by definition of the matrix $A_{i}$
\[
r(A_{i}x)=rH_{i}(\varphi),
\]
where
\[
H_{i}(\varphi)=\sqrt{\left(a^{(i)}_{11}\cos\varphi+a^{(i)}_{12}\sin\varphi\right)^{2}+
\left(a^{(i)}_{21}\cos\varphi+a^{(i)}_{22}\sin\varphi\right)^{2}}.
\]

Similarly, by definition of the matrix $A_{i}$
\[
\varphi(A_{i}x)=\Phi_{i}(\varphi),
\]
where
\[
\Phi_{i}(\varphi)=\arctan\left(\frac{a^{(i)}_{21}\cos\varphi+a^{(i)}_{22}\sin\varphi}%
{a^{(i)}_{11}\cos\varphi+a^{(i)}_{12}\sin\varphi}%
\right).
\]

From the obtained relations it follows that the first two
equalities in (\ref{E-lohibounds}) take the form
\[
\rho^{+}_{n}=\max_{\varphi}
\max_{i}\frac{H_{i}(\varphi)R_{n}(\Phi_{i}(\varphi))}{R_{n}(\varphi)},\quad
\rho^{-}_{n}=\min_{\varphi}
\max_{i}\frac{H_{i}(\varphi)R_{n}(\Phi_{i}(\varphi))}{R_{n}(\varphi)},\quad
\]
or, what is the same,
\begin{equation}\label{E-Rlohibounds}
\rho^{+}_{n}=\max_{\varphi}
\frac{R^{*}_{n}(\varphi)}{R_{n}(\varphi)},\quad
\rho^{-}_{n}=\min_{\varphi}
\frac{R^{*}_{n}(\varphi)}{R_{n}(\varphi)},
\end{equation}
where
\begin{equation}\label{E-Rstar}
R^{*}_{n}(\varphi)=\max_{i}\left\{H_{i}(\varphi)R_{n}(\Phi_{i}(\varphi))\right\}.
\end{equation}
The relations (\ref{E-auxnorm}) take the form
\[
rR_{n+1}(\varphi)=
 \max\left\{rR_{n}(\varphi),
 ~r\gamma^{-1}_{n}R^{*}_{n}(\varphi)\right\}
\]
or, equivalently,
\begin{equation}\label{E-Rauxnorm}
R_{n+1}(\varphi)=
 \max\left\{R_{n}(\varphi),
 ~\gamma^{-1}_{n}R^{*}_{n}(\varphi)\right\}
\end{equation}
and the normalization condition (\ref{E-newnorm}) takes the
form
\[
rR^{\circ}_{n+1}(\varphi)=\frac{rR_{n+1}(\varphi)}{r_{e}R_{n+1}(\varphi_{e})},
\]
where $(r_{e},\varphi_{e})$ are polar coordinates of the vector  $e$. Taking in place of $e$ the vector with polar coordinates $(1,0)$ the normalization condition can be rewritten in the form
\begin{equation}\label{E-Rnewnorm}
    R^{\circ}_{n+1}(\varphi)=\frac{R_{n+1}(\varphi)}{R_{n+1}(0)}.
\end{equation}

So, the max-relaxation iteration scheme can be represented as
follows. Given an averaging function $\gamma(\cdot,\cdot)$, set
$R_{0}(\varphi)\equiv 1$, and build recursively the
$2\pi$--periodic functions $R_{n}(\varphi)$ and
$R^{\circ}_{n}(\varphi)$, $n=1,2,\ldots$, in accordance with
the following rules:
\medskip

{\sl\begin{description}
\item[MR1] regarding the function $R_{n}(\varphi)$ already
    known compute the numerical values $\rho^{+}_{n}$ and
    $\rho^{-}_{n}$ by formulas (\ref{E-Rlohibounds}),
    (\ref{E-Rstar}) and set
    $\gamma_{n}=\gamma(\rho^{-}_{n},\rho^{+}_{n})$;

\item[MR2] define $R_{n+1}(\varphi)$ by (\ref{E-Rauxnorm}) and
    $R^{\circ}_{n+1}(\varphi)$ by (\ref{E-Rnewnorm}), and then
    determine the norm $\|\cdot\|^{\circ}_{n+1}$ as
    $\|x\|^{\circ}_{n+1}=rR^{\circ}_{n+1}(\varphi)$, where
    $(r,\varphi)$ are the polar coordinates of the vector $x$.
\end{description}}
\medskip

\begin{figure}[!htbp]
\begin{center}
\hfill\includegraphics[height=0.45\textwidth,clip]{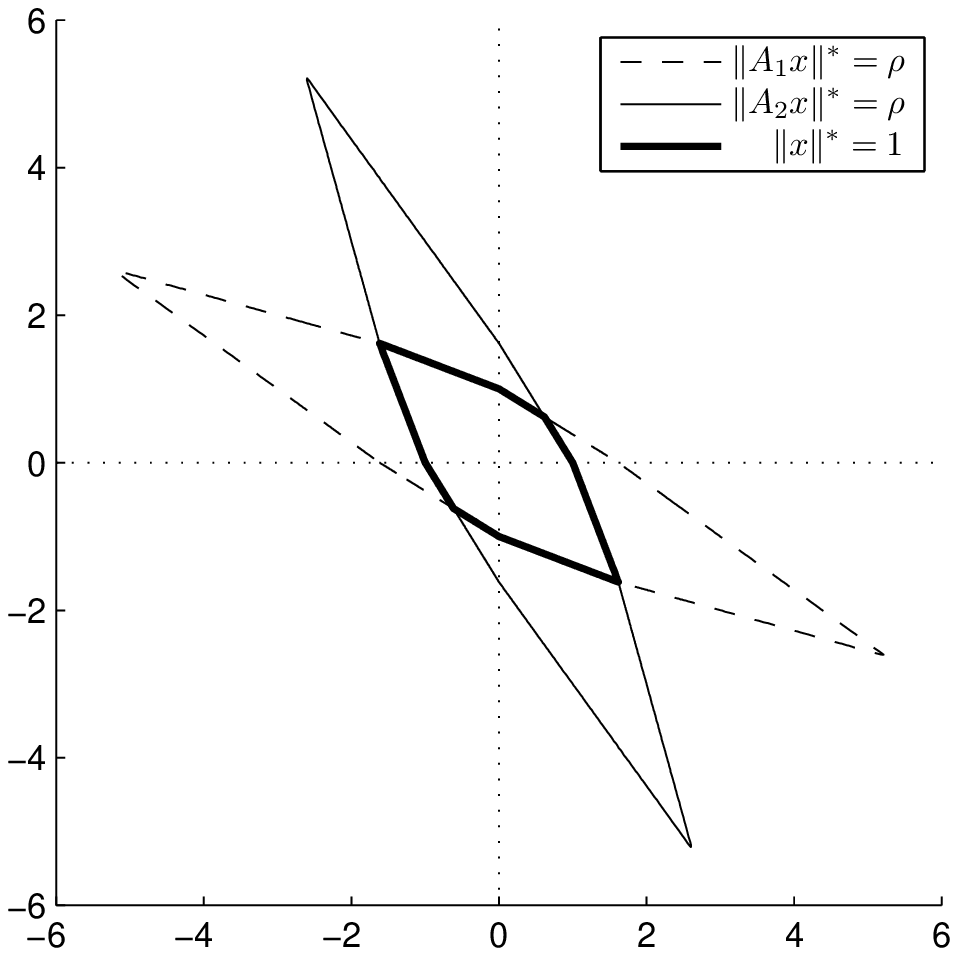}\hfill
\includegraphics[height=0.45\textwidth,clip]{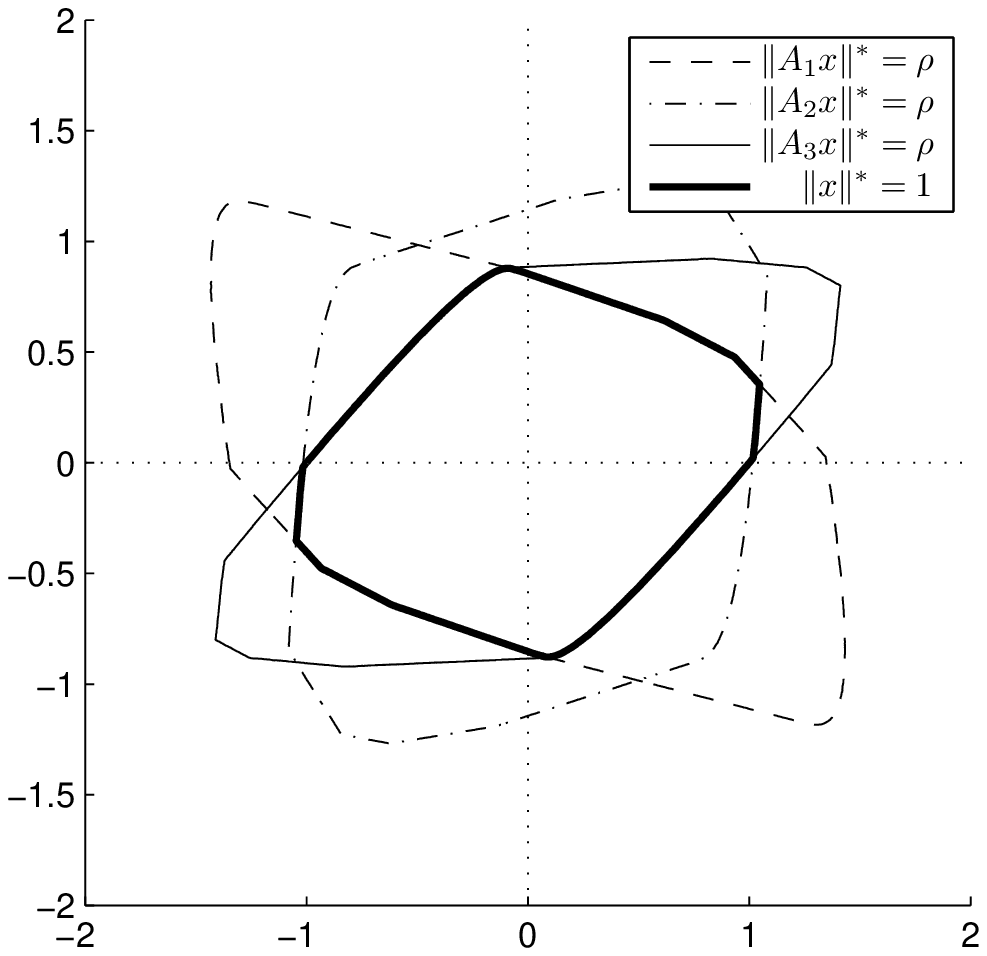}\hfill~
\caption{Examples of computation of Barabanov norms for $2\times 2$ matrices.}\label{F-barnorms}
\end{center}
\end{figure}

\begin{example}\label{Ex1}\rm
Consider the family $\setA=\{A_{1},A_{2}\}$ of $2\times 2$ matrices
\[
A_{1}=\left(\begin{array}{rr}
1&~1\\0&~1
\end{array}\right),\quad
A_{2}=\left(\begin{array}{rr}
1&~0\\1&~1
\end{array}\right).
\]
The functions $\Phi_{i}(\varphi),
H_{i}(\varphi),R_{n}(\varphi), R^{*}_{n}(\varphi)$ were chosen
to be piecewise linear with $3000$ nodes uniformly distributed
over the interval $[-\pi,\pi]$. It was needed $13$ iterations
of algorithm \textbf{MR1-MR2} with the averaging function
$\gamma(t,s)=\frac{t+s}{2}$ implemented in MATLAB to compute
the joint spectral radius $\rho(\setA)$ with the absolute
accuracy $10^{-3}$. The computed value of the joint spectral
radius is $\rho(\setA)=1.617$. The computed unit sphere of the
Barabanov norm $\|\cdot\|^{*}$ after the $13$th iteration of
algorithm \textbf{MR1-MR2} is shown on Fig.~\ref{F-barnorms} on
the left.
\end{example}

As is seen from Fig.~\ref{F-barnorms}, in Example~\ref{Ex1} the
sets $\|A_{1}x\|=\rho$ and $\|A_{2}x\|=\rho$ have exactly $4$
intersection points. This was theoretically proved in
\cite{Koz:CDC05:e,Koz:INFOPROC06:e} for the case when one of
the matrices $A_{1},A_{2}$ is lower triangle and the other is
upper triangle, and their entries are nonnegative. In
\cite{Koz:CDC05:e,Koz:INFOPROC06:e} this fact was one of key
points in disproving the Finiteness Conjecture. We do not know
whether this fact is valid in a general case or not, but
numerical tests based on algorithm \textbf{MR1-MR2} with
several dozens pairs of matrices $A_{1},A_{2}$ testify for this
fact.

\begin{example}\label{Ex2}\rm
Consider the family $\setA=\{A_{1},A_{2}, A_{3}\}$ of $2\times
2$ matrices
\[
A_{1}=\left(\begin{array}{rr}
1&~1\\0&~1
\end{array}\right),\quad
A_{2}=\left(\begin{array}{rr}
0.8&~0.6\\-0.6&~0.8
\end{array}\right),\quad
A_{3}=\left(\begin{array}{rr}
1\hphantom{.0}&0\hphantom{.0}\\-0.4&1.3
\end{array}\right).
\]
Here the functions $\Phi_{i}(\varphi),
H_{i}(\varphi),R_{n}(\varphi), R^{*}_{n}(\varphi)$ were also
chosen to be piecewise linear with $3000$ nodes uniformly
distributed over the interval $[-\pi,\pi]$. It was needed $31$
iterations of algorithm \textbf{MR1-MR2} with the averaging
function $\gamma(t,s)=\frac{t+s}{2}$ implemented in MATLAB to
compute the joint spectral radius $\rho(\setA)$ with the
absolute accuracy $10^{-3}$. The computed value of the joint
spectral radius is $\rho(\setA)=1.347$. The computed unit
sphere of the Barabanov norm $\|\cdot\|^{*}$ after the $31$d
iteration of algorithm \textbf{MR1-MR2} is shown on
Fig.~\ref{F-barnorms} on the right. The MATLAB code for this
Example can be found in \cite{Koz:ArXiv10-1}.
\end{example}

\section{Concluding remarks}\label{S-rem}

The max-relaxation algorithm suggested in the paper allows to
calculate the joint spectral radius of a finite matrix family
(of arbitrary matrix size and arbitrary amount of matrices in
the set) with any required accuracy and to evaluate a
posteriori the computational error. Still, this algorithm gives
rise to a set of open problems.

\begin{problem}While the quantities
$\{\rho^{\pm}_{n}\}$ provide a posteriori bounds for the
accuracy of approximation of $\rho({\setA})$ the question about
the accuracy of approximation of the Barabanov norm
$\|\cdot\|^{*}$ by the norms $\|\cdot\|^{\circ}_{n}$ is open.
\end{problem}

It seems, the difficulty in resolving this problem is caused by
the fact that in general the Barabanov norms for a matrix
family are determined ambiguously. Namely to overcome this
difficulty we preferred to consider relaxation algorithm
instead of direct one. Moreover, as can be shown, both
theoretically and  numerically, if to set $\|x\|_{n+1}=
\gamma^{-1}_{n}\max_{i}\|A_{i}x\|_{n}$ in (\ref{E-auxnorm})
then the obtained direct computational analog of algorithm
\textbf{MR1-MR2} may turn out to be non-convergent.

\begin{problem}
The question about the rate of convergence of the sequences
$\{\rho^{+}_{n}\}$ and $\{\rho^{-}_{n}\}$ to the joint spectral
radius is also open.
\end{problem}

Remark also that in this paper mainly the algorithm for
building of Barabanov norms rather than its computational
details was studied. The numerical aspects of implementation of
this algorithm require additional analysis.

\begin{problem}
An estimation of the computational cost of the max-relaxation
algorithm is required. The dependance of the algorithm on
parameters $r$, $m$ and the choice of the averaging function
$\gamma(t,s)$ is acute, too.
\end{problem}

\section*{Acknowledgements}
This work was supported by the Russian Foundation for Basic
Research, project no. 10-01-00175.

\end{document}